\documentclass{amsproc} 

\newcommand{\ttr}{\ensuremath{\triangleright}}

\usepackage{latexsym}
\usepackage{amsfonts}
\usepackage{amssymb}
\usepackage{amsmath}
\usepackage{amscd}
\usepackage{eucal}
\usepackage[all, knot]{xy}
\xyoption{arc}
\usepackage{amsthm}
\usepackage[dvips]{pict2e}

\usepackage{verbatim} 
\usepackage[dvips]{graphicx}

\newcommand{\w}{\omega}

\newtheorem{Def}{Definition}[section]
\newtheorem{thm}[Def]{Theorem}

\newtheorem{prop}[Def]{Proposition}
\newtheorem{cor}[Def]{Corollary}

\newtheoremstyle{example}{\topsep}{\topsep}%
     {}
     {}
     {\bfseries}
     {.}
     {2pt}
     {\thmname{#1}\thmnumber{ #2}\thmnote{ #3}}

   \theoremstyle{example}

\newtheorem{hyp}[Def]{Hypothesis}

\newcommand{\lra}{\ensuremath{\longrightarrow}}
\newcommand{\cat}[1]{\ensuremath{\mbox{\bfseries {\upshape {#1}}}}}
\newcommand{\cl}[1]{\ensuremath{\mathcal {#1}}}

\newcommand{\ed}{\end{document}}
\newcommand{\map}[1]{\ensuremath{\stackrel{{#1}}{\lra}}}

\newcommand{\Cat}{{\cat{Cat}}}

\newcommand{\numarabic}{\renewcommand{\labelenumi}{\arabic{enumi}.}}

\newcommand{\sunit}{\setlength{\unitlength}{1mm}}
\newcommand{\setunit}[1]{\setlength{\unitlength}{#1}}

\newenvironment{prf}{\vspace{2ex}\begin{sloppypar}{\noindent\upshape
{\bfseries Proof. }}} {{\hspace*{\fill}
$\Box$}\end{sloppypar}\vspace{2ex}}

\newcommand{\triarrow}{
\sunit
\begin{picture}(24,0)(-0.5,-1)
\put(0,0){\line(1,0){20}}
\put(19.8,-0.2){\makebox(0,0)[l]{\ttr}}
\end{picture}}

\newcommand{\degenarrow}{
\begin{picture}(0,0)
\qbezier[20](0,0)(-3,-2.25)(-6,-4.5)
\put(-6,-4.5){\vector(-12,-9){0.1}}
\end{picture}}

\newcommand{\rrr}{
\begin{picture}(60,70)(0,15)

\qbezier(2,9)(26,32)(35,60)
\qbezier(35,60)(36,62)(35.6,64)
\qbezier(35.6,64)(35.2,68.2)(32.2,65)
\qbezier(32.2,65)(29,61.5)(28.5,55)
\qbezier(28.5,55)(27.4,42.4)(41,46)
\qbezier(41,46)(50.2,48.8)(58.4,57)
\qbezier(58.4,57)(51.6,40)(52,28)
\qbezier(52,28)(52.6,12.2)(64.8,9.5)
\qbezier(64.8,9.5)(71.5,7.8)(77,14)
\end{picture}
}

\newcommand{\curlyr}{\setunit{0.03mm}
\makebox[2.8mm][l]{\begin{picture}(110,70)(38,-10)
\put(0,0){\rrr}
\put(1,0){\rrr}
\put(2,0){\rrr}
\end{picture}}}

\newcommand{\smallcurlyr}{\setunit{0.027mm}
\makebox[2.7mm][l]{\begin{picture}(110,70)(38,-10)
\put(0,0){\rrr}
\put(1,0){\rrr}
\put(2,0){\rrr}
\end{picture}}}

\newcommand{\aaa}{
\begin{picture}(80,50)
\qbezier(2,9)(3.8,5.5)(10,5)
\qbezier(10,5)(24.5,4.8)(40,12.8)
\qbezier(40,12.8)(47.5,16.5)(50.8,20)
\qbezier(50.8,20)(56.5,26)(56,36)
\qbezier(56,36)(55,50)(41,45)
\qbezier(41,45)(31.5,42)(25,31)
\qbezier(25,31)(19,21)(25,10)
\qbezier(25,10)(30,1)(41,8)
\qbezier(41,8)(46.5,11.2)(49.5,17.5)
\qbezier(49.5,17.5)(52,6.8)(61,5)
\qbezier(61,5)(69.2,3.5)(75,9)
\end{picture}}

\newcommand{\curlya}{\setunit{0.04mm}
\makebox[3.5mm][l]{\begin{picture}(80,70)(28,5)
\put(0,0){\aaa}
\put(1,0){\aaa}
\put(2,0){\aaa}
\end{picture}}}

\newcommand{\smallcurlya}{\setunit{0.036mm}
\makebox[3.3mm][l]{\begin{picture}(80,70)(28,5)
\put(0,0){\aaa}
\put(1,0){\aaa}
\put(2,0){\aaa}
\end{picture}}}

\newcommand{\llll}{
\begin{picture}(65,80)
\qbezier(5,7)(9,4)(15,6)
\qbezier(15,6)(23.5,8.5)(29,15)
\qbezier(29,15)(38,26)(42.8,40)
\qbezier(42.8,40)(46,49)(47,60)
\qbezier(47,60)(48.5,65.6)(47.8,71)
\qbezier(47.8,71)(47.5,77.5)(44,72)
\qbezier(44,72)(41.5,68)(40.5,60)
\qbezier(40.5,60)(37.5,41)(38.5,30)
\qbezier(38.5,30)(38.5,22)(40,18)
\qbezier(40,18)(42,12)(46,8)
\qbezier(46,8)(52,1.5)(58.5,8)
\end{picture}}

\newcommand{\curlyl}{\setunit{0.04mm}
\makebox[2.9mm][l]{\begin{picture}(80,70)(28,5)
\put(0,0){\llll}
\put(1,0){\llll}
\put(2,0){\llll}
\end{picture}}}

\newcommand{\smallcurlyl}{\setunit{0.036mm}
\makebox[2.8mm][l]{\begin{picture}(80,70)(28,5)
\put(0,0){\llll}
\put(1,0){\llll}
\put(2,0){\llll}
\end{picture}}}

\newcommand{\catk}[2]{\ensuremath{\cat{{#1}}({\ensuremath{#2}})}} 
\newcommand{\catkj}[3]{\ensuremath{\catk{#1}{#2}_{#3}}}

\newcommand{\ncat}{$n$-category}
\newcommand{\ncats}{$n$-categories}

\newcommand{\pdfshift}{\setlength{\topmargin}{25mm}}

\pdfshift

\begin{document}

\title[Degenerate categories and degenerate bicategories]{The periodic table of $n$-categories for low dimensions I: degenerate categories and degenerate bicategories}
\author{Eugenia Cheng}
\address{Department of Mathematics, University of Chicago, Chicago, IL60637}
\email{eugenia@math.uchicago.edu}
\author{Nick Gurski}
\address{Department of Mathematics, University of Chicago, Chicago, IL60637}
\email{gurski@math.uchicago.edu}


\date{July 2005}

\begin{abstract}

We examine the periodic table of weak $n$-categories for the low-dimensional cases.  It is widely understood that degenerate categories give rise to monoids, doubly degenerate bicategories to commutative monoids, and degenerate bicategories to monoidal categories; however, to understand this correspondence fully we examine the totalities of such structures together with maps between them and higher maps between those. Categories naturally form a 2-category {\bfseries Cat} so we take the full sub-2-category of this whose 0-cells are the degenerate categories.  Monoids naturally form a category, but we regard this as a discrete 2-category to make the comparison.  We show that this construction does not yield a biequivalence; to get an equivalence we ignore the natural transformations and consider only the {\it category} of degenerate categories. A similar situation occurs for degenerate bicategories.  The tricategory of such does not yield an equivalence with monoidal categories; we must consider only the categories of such structures.
For doubly degenerate bicategories the tricategory of such is not naturally triequivalent to the category of commutative monoids (regarded as a tricategory).  However in this case considering just the categories does not give an equivalence either; to get an equivalence we consider the {\it bicategory} of doubly degenerate bicategories. We conclude with a hypothesis about how the above cases might generalise for $n$-fold degenerate $n$-categories.

\end{abstract}

\maketitle

\section*{Introduction}


In this paper we examine the first few entries in the ``Periodic Table" of $n$-categories.  This table was first described by Baez and Dolan in \cite{bd3} and is closely linked to the Stabilisation Hypothesis.

The idea of the Periodic Table is to study ``degenerate" forms of $n$-category, that is, $n$-categories that are trivial below a certain dimension $k$.  Now, such an $n$-category only has non-trivial cells in the top $n-k$ dimensions, so we can perform a ``dimension shift" and regard this as an $(n-k)$-category: the old $k$-cells become the new 0-cells, the old $(k+1)$-cells become the new 1-cells, and so on up to the old $n$-cells which become the new $(n-k)$-cells.  We call this a ``$k$-fold degenerate $n$-category", and the dimension-shift is depicted in the schematic diagram in Figure~\ref{dshift}.

\vspace{5cm}

\begin{figure}
\caption{Dimension-shift for $k$-fold degenerate $n$-categories} \label{dshift}
\setunit{1mm}
\begin{tabular}{|lcc|} \hline &&\\[-4pt]
{\bfseries ``old" $n$-category} & \triarrow & {\bfseries ``new" $(n-k)$-category} \\[8pt] \hline && \\[-4pt]
\hspace{1em}$\left.\begin{array}{l}
\mbox{0-cells} \\ \mbox{1-cells} \\ \hspace{1em}\vdots \\[4pt] \makebox(13,0)[r]{$(k-1)$-cells}
\end{array}\right\}$ trivial &&\\[32pt]
\hspace{1.4em}$k$-cells & \triarrow & 0-cells \\[6pt]
$(k+1)$-cells & \triarrow & 1-cells \\
\hspace{2.2em} \vdots & \vdots & \vdots \\[2pt]
\hspace{1.4em}$n$-cells & \triarrow & $(n-k)$-cells \\[4pt] \hline
\end{tabular}
\end{figure}

However, this process evidently yields a special kind of $(n-k)$-category.  This is because the 0-cells in the ``new" $(n-k)$-category have some extra structure on them, which comes from all the different types of composition they had as $k$-cells in the ``old" $n$-category.  Essentially, we get one type of ``multiplication" (or tensor) for each type of composition that there was, and these different tensors interact according to the old interchange laws for composition.  Since $k$-cells have $k$ types of composition, we have $k$ different monoidal structures; this is what is known as a ``$k$-tuply monoidal $(n-k)$-category" although the precise general definition has not been made.

\begin{figure}
\caption{The hypothesised Periodic Table of $n$-categories} \label{ptable}

\setunit{2mm}  \nopagebreak[4]
\begin{picture}(80,85)(16,-5)

\put(0,37.5){
\renewcommand{\tabcolsep}{1.5em}

\begin{tabular}{|lllll|} \hline &&&&\\
set & category & 2-category & 3-category & \makebox(0,0)[r]{$\cdots$}  \\[40pt]
monoid & monoidal category & monoidal 2-category & monoidal 3-category & \makebox(0,0)[r]{$\cdots$} \\
\makebox(2.5,0.5)[l]{$\equiv$}category with & \makebox(2.5,0.5)[l]{$\equiv$}2-category with & \makebox(2.5,0.5)[l]{$\equiv$}3-category with & \makebox(2.5,0.5)[l]{$\equiv$}4-category with &\\ 
\makebox(2.5,0.5)[l]{}only one object & \makebox(2.5,0.5)[l]{}only one object & \makebox(2.5,0.5)[l]{}only one object & \makebox(2.5,0.5)[l]{}only one object &\\[35pt]
{\bfseries commutative}& braided monoidal & braided monoidal & braided monoidal & \makebox(0,0)[r]{$\cdots$}  \\
{\bfseries monoid}& category & 2-category & 3-category &\\
\makebox(2.5,0.5)[l]{$\equiv$}2-category with & \makebox(2.5,0.5)[l]{$\equiv$}3-category with & \makebox(2.5,0.5)[l]{$\equiv$}4-category with & \makebox(2.5,0.5)[l]{$\equiv$}5-category with &\\ 
\makebox(2.5,0.5)[l]{}only one object & \makebox(2.5,0.5)[l]{}only one object & \makebox(2.5,0.5)[l]{}only one object & \makebox(2.5,0.5)[l]{}only one object &\\
\makebox(2.5,0.5)[l]{}only one 1-cell & \makebox(2.5,0.5)[l]{}only one 1-cell & \makebox(2.5,0.5)[l]{}only one 1-cell & \makebox(2.5,0.5)[l]{}only one 1-cell &\\[30pt]
\makebox(8,0)[r]{\large $\prime\prime$}& \makebox(16,0)[l]{{\bfseries symmetric monoidal}} & sylleptic monoidal & \makebox(8,0)[r]{\Large ?} & \makebox(0,0)[r]{$\cdots$}  \\
& {\bfseries category} & 2-category & &\\
\makebox(2.5,0.5)[l]{$\equiv$}3-category with & \makebox(2.5,0.5)[l]{$\equiv$}4-category with & \makebox(2.5,0.5)[l]{$\equiv$}5-category with & \makebox(2.5,0.5)[l]{$\equiv$}6-category with &\\ 
\makebox(2.5,0.5)[l]{}only one object & \makebox(2.5,0.5)[l]{}only one object & \makebox(2.5,0.5)[l]{}only one object & \makebox(2.5,0.5)[l]{}only one object &\\
\makebox(2.5,0.5)[l]{}only one 1-cell & \makebox(2.5,0.5)[l]{}only one 1-cell & \makebox(2.5,0.5)[l]{}only one 1-cell & \makebox(2.5,0.5)[l]{}only one 1-cell &\\
\makebox(2.5,0.5)[l]{}only one 2-cell & \makebox(2.5,0.5)[l]{}only one 2-cell & \makebox(2.5,0.5)[l]{}only one 2-cell & \makebox(2.5,0.5)[l]{}only one 2-cell &\\ [20pt]
\makebox(8,0)[r]{\large $\prime\prime$}& \makebox(8,0)[r]{\large $\prime\prime$} & \makebox(16,0)[l]{{\bfseries symmetric monoidal}} & \makebox(8,0)[r]{\Large ?}& \makebox(0,0)[r]{$\cdots$}  \\
& & {\bfseries 2-category} & &\\
\makebox(2.5,0.5)[l]{$\equiv$}4-category with & \makebox(2.5,0.5)[l]{$\equiv$}5-category with & \makebox(2.5,0.5)[l]{$\equiv$}6-category with & \makebox(2.5,0.5)[l]{$\equiv$}7-category with &\\ 
\makebox(2.5,0.5)[l]{}only one object & \makebox(2.5,0.5)[l]{}only one object & \makebox(2.5,0.5)[l]{}only one object & \makebox(2.5,0.5)[l]{}only one object &\\
\makebox(2.5,0.5)[l]{}only one 1-cell & \makebox(2.5,0.5)[l]{}only one 1-cell & \makebox(2.5,0.5)[l]{}only one 1-cell & \makebox(2.5,0.5)[l]{}only one 1-cell &\\
\makebox(2.5,0.5)[l]{}only one 2-cell & \makebox(2.5,0.5)[l]{}only one 2-cell & \makebox(2.5,0.5)[l]{}only one 2-cell & \makebox(2.5,0.5)[l]{}only one 2-cell &\\ 
\makebox(2.5,0.5)[l]{}only one 3-cell & \makebox(2.5,0.5)[l]{}only one 3-cell & \makebox(2.5,0.5)[l]{}only one 3-cell & \makebox(2.5,0.5)[l]{}only one 3-cell &\\[30pt] 
\makebox(8,0)[r]{\large $\prime\prime$}& \makebox(8,0)[r]{\large $\prime\prime$}& \makebox(8,0)[r]{\large $\prime\prime$}& \makebox(8,0)[r]{\Large ?}& \makebox(0,0)[r]{$\cdots$} \\
\makebox(7.5,0)[r]{$\vdots$} & \makebox(7.5,0)[r]{$\vdots$}
& \makebox(7.5,0)[r]{$\vdots$} & \makebox(7.5,0)[r]{$\vdots$} & \\
&&&& \\
&&&& \\
 
\hline
\end{tabular}}

\put(22,76){\degenarrow}
\put(22,61){\degenarrow}
\put(22,46){\degenarrow}
\put(22,30){\degenarrow}
\put(22,12){\degenarrow}

\put(43.5,76){\degenarrow}
\put(43.5,61){\degenarrow}
\put(43.5,46){\degenarrow}
\put(43.5,30){\degenarrow}
\put(43.5,12){\degenarrow}

\put(64,76){\degenarrow}
\put(64,61){\degenarrow}
\put(64,46){\degenarrow}
\put(64,30){\degenarrow}
\put(64,12){\degenarrow}


\end{picture}
\end{figure}
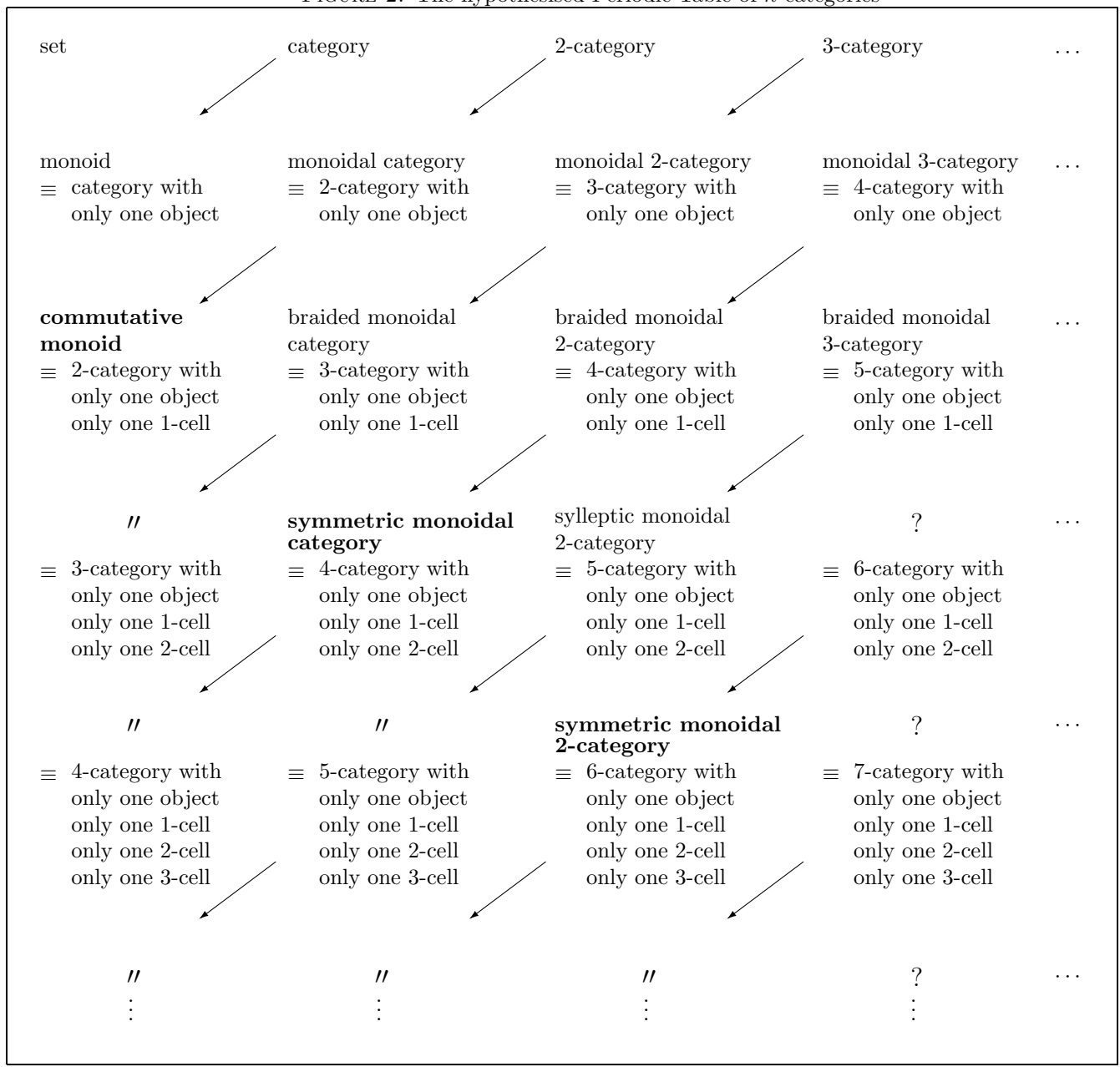

A natural question to ask then is: exactly what sort of $(n-k)$-category structure does this degeneracy process produce?  This is the question that the Periodic Table seeks to answer.  Figure~\ref{ptable} shows the first few columns of the hypothesised Periodic Table.  

(In this table we follow Baez and Dolan and omit the word ``weak" understanding that all the $n$-categories in consideration are weak.)
The dotted arrows indicate the process of collapsing the lowest dimension of the structure, that is, considering the one-object case. The entries in bold face show where the table is supposed to have ``stabilised" --- the entries in the rest of the column underneath should continue to be the same, according to the Stabilisation Hypothesis \cite{bd3}.

Thus the first column says: 

\begin{itemize}
\item A category with only one object ``is" a monoid.
\item A bicategory with only one 0-cell and only one 1-cell ``is" a commutative monoid.
\item A tricategory with only one 0-cell, 1-cell and 2-cell ``is" a commutative monoid.
\end{itemize}

\enlargethispage{1em}

The second column says:
\begin{itemize}
\item A bicategory with only one 0-cell ``is" a monoidal category.
\item A tricategory with only one 0-cell and 1-cell ``is" a braided monoidal category.
\item A tetracategory with only one 0-cell, 1-cell and 2-cell ``is" a symmetric monoidal category.
\end{itemize}

In each case we need to say exactly what ``is" means.   
In this paper we examine the top left hand corner of the Periodic Table, that is, degenerate categories and bicategories.  (In a future paper we will examine tricategories.)    

The main problem is the presence of some unwanted extra structure in the ``new" $(n-k)$-categories in the form of distinguished elements, arising from the structure constraints in the original $n$-categories --- a specified $k$-cell structure constraint in the ``old" $n$-category will appear as a distinguished 0-cell in the ``new" $(n-k)$-category under the dimension-shift depicted in Figure~\ref{dshift}.   We will show that some care is thus required in the interpretion of the above statements. (For $n=2$ this phenomenon is mentioned by Leinster in \cite{lei4} and was further described in a talk \cite{lei10}.) 

We begin in Section~\ref{method} by outlining our methodology.  In Section~\ref{dcs} we describe the well-known example of degenerate categories; in this case the 1-cells form a monoid with multiplication given by composition.

In Section~\ref{vdbs} we examine ``doubly degenerate" bicategories, that is, bicategories with only one 0-cell and 1-cell.  Now the 2-cells have two compositions on them---horizontal and vertical.  So we might expect the 2-cells to form some sort of structure with two different multiplications; however, we can use an Eckmann-Hilton argument to show that these two multiplications are the same and in fact commutative.  In Section~\ref{dbs} we study degenerate bicategories.  Here, the 1-cells become the objects of a monoidal category, with tensor given by the old composition of 1-cells.  

These basic results are to some extent well-known \cite{gps1, lei8}, but the focus of this paper is to make a precise interpretation of the statements in the Periodic Table by examining the \emph{totality} of each of the above structures, in the sense that we discuss in Section~\ref{method}.  We sum up the results as follows.

\enlargethispage{5cm}

\begin{figure}[b]
\caption{Comparison of overall structure for degenerate categories}
\label{sdone}
\setunit{2mm}  \nopagebreak[4]
\begin{picture}(80,38)(11,-2)

\put(2,0){\framebox(78,35){}}

\put(5,15){\makebox(0,0)[l]{degenerate categories}}
\put(5,10){\makebox(0,0)[l]{functors}}
\put(5,5){\makebox(0,0)[l]{natural transformations}}

\put(22,12){\makebox(0,0)[l]{$\left.\begin{array}{l}
\\[14mm] \end{array}\right\}$ \textsf{category}}}

\put(31,10){\makebox(0,0)[l]{$\left.\begin{array}{l}
\\[22mm] \end{array}\right\}$ \textsf{bicategory}}}

\put(65,15){\makebox(0,0)[l]{monoids}}
\put(65,10){\makebox(0,0)[l]{homomorphisms}}
\put(65,5){\makebox(0,0)[l]{identities}}

\put(66.5,12){\makebox(0,0)[r]{\textsf{category} $\left\{\begin{array}{l}
\\[14mm] \end{array}\right.$}}

\put(58,10){\makebox(0,0)[r]{\textsf{bicategory} $\left\{\begin{array}{l}
\\[22mm] \end{array}\right.$}}

\qbezier[15](39,11.5)(41,15)(44.5,15)
\qbezier[15](44.5,15)(48,15)(50,11.5)
\put(50,11.5){\vector(2,-3){0.5}}

\put(44.5,16){\makebox(0,0)[c]{\bfseries\textsf{no obvious functor}}}

\qbezier(29,13.5)(34,22)(44,22)
\qbezier(44,22)(54,22)(59,13.5)
\put(59,13.5){\vector(2,-3){0.5}}
\put(44,23){\makebox(0,0)[c]{\bfseries\textsf{equivalence}}}

\put(5,32){\makebox(0,0)[l]{\sc{Totality of}}}
\put(5,29.5){\makebox(0,0)[l]{\sc{degenerate categories}}}

\put(65,32){\makebox(0,0)[l]{\sc{Totality of}}}
\put(65,29.5){\makebox(0,0)[l]{\sc{monoids}}}

\put(32,31){\vector(1,0){26}}
\put(45,32){\makebox(0,0)[c]{\textsf{structure comparison}}}

\end{picture}
\end{figure}
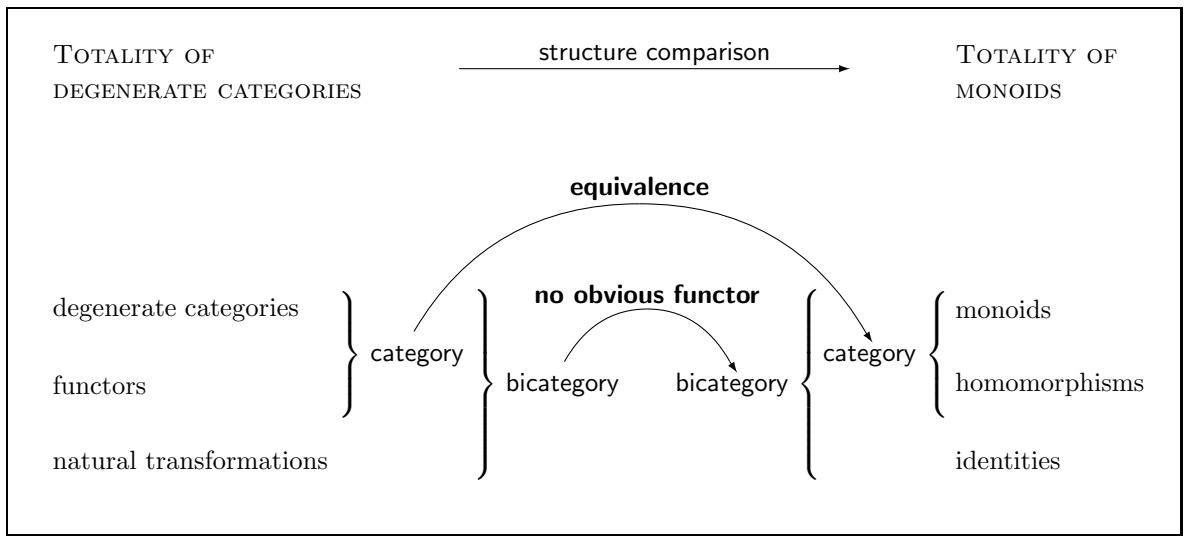

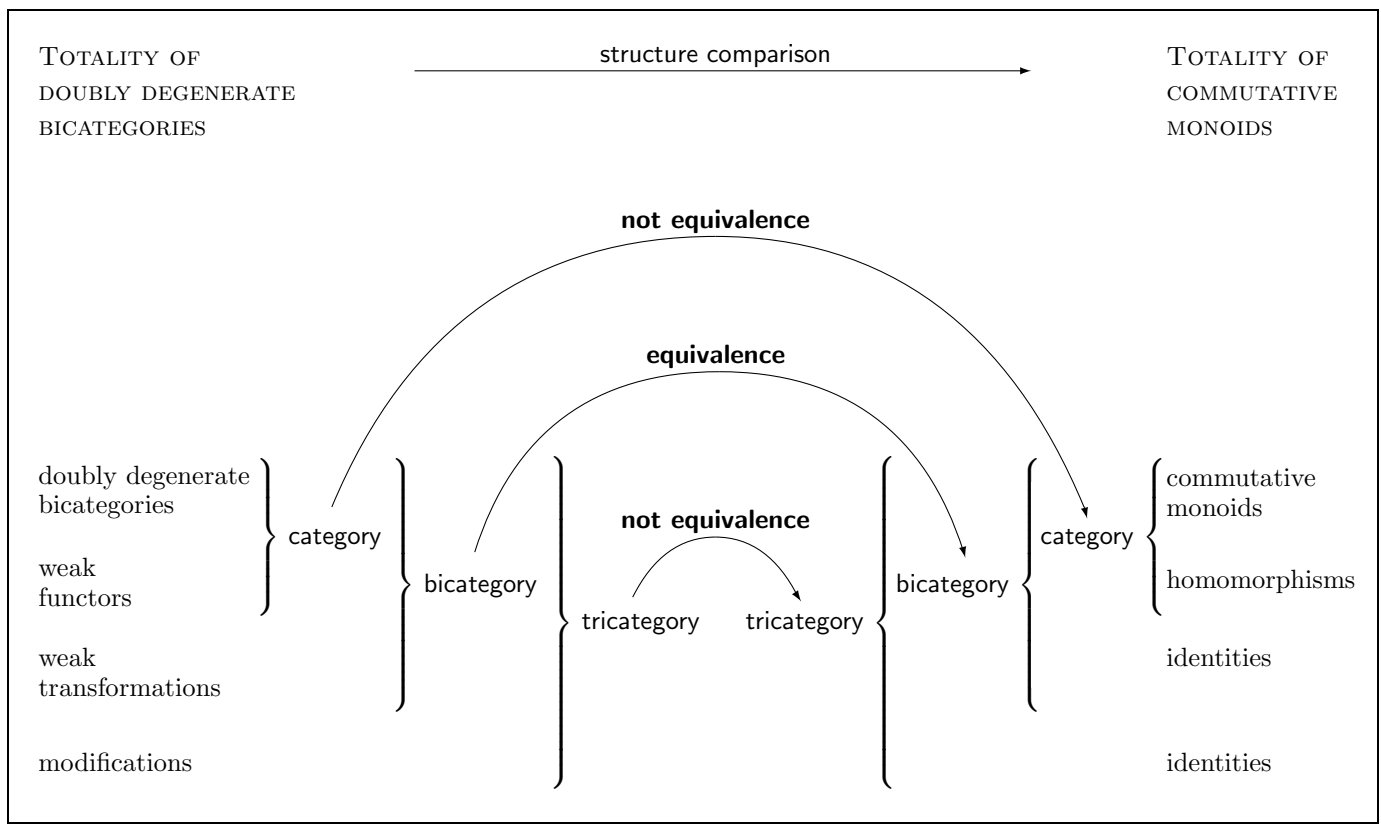
\begin{figure}
\caption{Comparison of overall structure for doubly degenerate bicategories}
\label{sdtwo}
\setunit{2mm}  
\begin{picture}(88,56)(14,0)

\put(-2,0){\framebox(91,54){}}

\put(0,23){\makebox(0,0)[l]{doubly degenerate}}
\put(0,21){\makebox(0,0)[l]{bicategories}}
\put(0,17){\makebox(0,0)[l]{weak}}
\put(0,15){\makebox(0,0)[l]{functors}}
\put(0,11){\makebox(0,0)[l]{weak}}
\put(0,9){\makebox(0,0)[l]{transformations}}
\put(0,4){\makebox(0,0)[l]{modifications}}

\put(12.5,19){\makebox(0,0)[l]{$\left.\begin{array}{l}
\\[17mm] \end{array}\right\}$ \textsf{category}}}

\put(21.5,15.8){\makebox(0,0)[l]{$\left.\begin{array}{l}
\\[29.5mm] \end{array}\right\}$ \textsf{bicategory}}}

\put(32,13.3){\makebox(0,0)[l]{$\left.\begin{array}{l}
\\[40mm] \end{array}\right\}$ \textsf{tricategory}}}

\put(75,23){\makebox(0,0)[l]{commutative}}
\put(75,21){\makebox(0,0)[l]{monoids}}
\put(75,16){\makebox(0,0)[l]{homomorphisms}}
\put(75,11){\makebox(0,0)[l]{identities}}
\put(75,4){\makebox(0,0)[l]{identities}}

\put(66.6,19){\makebox(0,0)[l]{\textsf{category} $\left\{\begin{array}{l}
\\[17mm] \end{array}\right.$}}

\put(57,15.8){\makebox(0,0)[l]{\textsf{bicategory} $\left\{\begin{array}{l}
\\[29.5mm] \end{array}\right.$}}

\put(47,13.3){\makebox(0,0)[l]{\textsf{tricategory} $\left\{\begin{array}{l}
\\[40mm] \end{array}\right.$}}

\qbezier(39.5,15)(41.5,19)(45,19)
\qbezier(45,19)(48.5,19)(50.5,15)
\put(50.5,15){\vector(2,-3){0.2}}

\qbezier(29,18)(32.8,30)(45.25,30)
\qbezier(45.25,30)(57.7,30)(61.5,18)
\put(61.5,18){\vector(1,-4){0.1}}

\qbezier(19.5,21)(27,39)(45,39)
\qbezier(45,39)(62,39)(69.5,21)
\put(69.5,21){\vector(1,-4){0.2}}

\put(45,20){\makebox(0,0)[c]{\bfseries\textsf{not equivalence}}}
\put(45,31){\makebox(0,0)[c]{\bfseries\textsf{equivalence}}}
\put(45,40){\makebox(0,0)[c]{\bfseries\textsf{not equivalence}}}

\put(0,51){\makebox(0,0)[l]{\sc{Totality of}}}
\put(0,48.5){\makebox(0,0)[l]{\sc{doubly degenerate}}}
\put(0,46.2){\makebox(0,0)[l]{\sc{bicategories}}}

\put(75,51){\makebox(0,0)[l]{\sc{Totality of}}}
\put(75,48.5){\makebox(0,0)[l]{\sc{commutative}}}
\put(75,46.2){\makebox(0,0)[l]{\sc{monoids}}}

\put(25,50){\vector(1,0){41}}
\put(45,51){\makebox(0,0)[c]{\textsf{structure comparison}}}

\end{picture}
\end{figure}

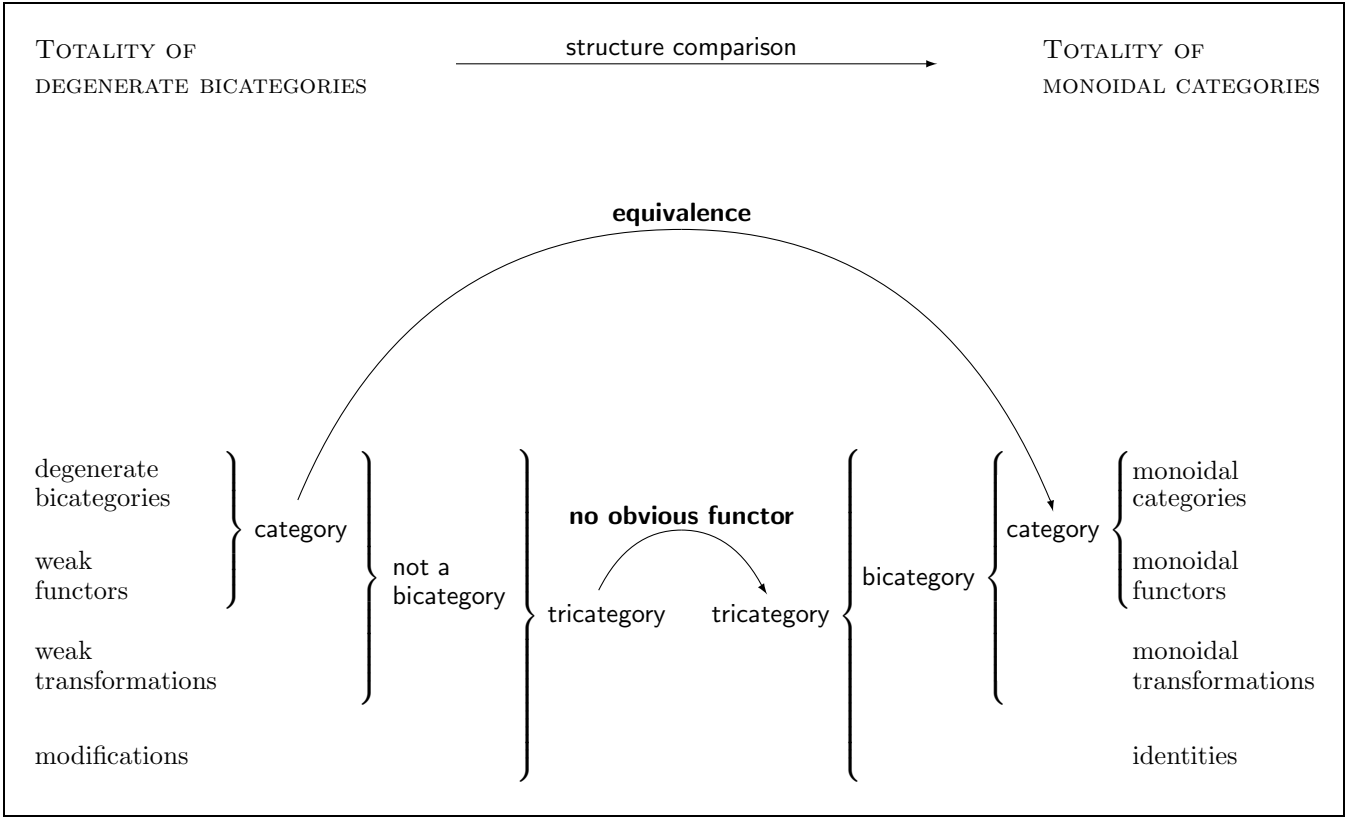
\begin{figure}
\caption{Comparison of overall structure for degenerate bicategories}
\label{sdthree}
\setunit{2mm}  \nopagebreak[4]
\begin{picture}(88,56)(14,0)

\put(0,0){\framebox(89,54){}}

\put(2,23){\makebox(0,0)[l]{degenerate}}
\put(2,21){\makebox(0,0)[l]{bicategories}}
\put(2,17){\makebox(0,0)[l]{weak}}
\put(2,15){\makebox(0,0)[l]{functors}}
\put(2,11){\makebox(0,0)[l]{weak}}
\put(2,9){\makebox(0,0)[l]{transformations}}
\put(2,4){\makebox(0,0)[l]{modifications}}

\put(12.5,19){\makebox(0,0)[l]{$\left.\begin{array}{l}
\\[17mm] \end{array}\right\}$ \textsf{category}}}

\put(21.5,15.8){\makebox(0,0)[l]{$\left.\begin{array}{l}
\\[29.5mm] \end{array}\right\}$}}

\put(25.8,16.5){\makebox(0,0)[l]{\textsf{not a}}}
\put(25.8,14.5){\makebox(0,0)[l]{\textsf{bicategory}}}

\put(32,13.3){\makebox(0,0)[l]{$\left.\begin{array}{l}
\\[40mm] \end{array}\right\}$ \textsf{tricategory}}}

\put(75,23){\makebox(0,0)[l]{monoidal}}
\put(75,21){\makebox(0,0)[l]{categories}}
\put(75,17){\makebox(0,0)[l]{monoidal}}
\put(75,15){\makebox(0,0)[l]{functors}}
\put(75,11){\makebox(0,0)[l]{monoidal}}
\put(75,9){\makebox(0,0)[l]{transformations}}
\put(75,4){\makebox(0,0)[l]{identities}}

\put(66.6,19){\makebox(0,0)[l]{\textsf{category} $\left\{\begin{array}{l}
\\[17mm] \end{array}\right.$}}

\put(57,15.8){\makebox(0,0)[l]{\textsf{bicategory} $\left\{\begin{array}{l}
\\[29.5mm] \end{array}\right.$}}

\put(47,13.3){\makebox(0,0)[l]{\textsf{tricategory} $\left\{\begin{array}{l}
\\[40mm] \end{array}\right.$}}

\qbezier[15](39.5,15)(41.5,19)(45,19)
\qbezier[15](45,19)(48.5,19)(50.5,15)
\put(50.5,15){\vector(2,-3){0.2}}


\qbezier(19.5,21)(27,39)(45,39)
\qbezier(45,39)(62,39)(69.5,21)
\put(69.5,21){\vector(1,-4){0.2}}

\put(45,40){\makebox(0,0)[c]{\bfseries\textsf{equivalence}}}
\put(45,20){\makebox(0,0)[c]{\bfseries\textsf{no obvious functor}}}

\put(2,51){\makebox(0,0)[l]{\sc{Totality of}}}
\put(2,48.5){\makebox(0,0)[l]{\sc{degenerate bicategories}}}

\put(69,51){\makebox(0,0)[l]{\sc{Totality of}}}
\put(69,48.5){\makebox(0,0)[l]{\sc{monoidal categories}}}

\put(30,50){\vector(1,0){32}}
\put(45,51){\makebox(0,0)[c]{\textsf{structure comparison}}}

\end{picture}
\end{figure}

\begin{itemize}

\item Comparing each degenerate category with the monoid formed by its 1-cells, we exhibit an equivalence of categories of these structures, but not a biequivalence of bicategories; see Figure~\ref{sdone}.

\item Comparing each doubly degenerate bicategory with the commutative monoid formed by its 2-cells, we exhibit a biequivalence of bicategories of these structures, but not an equivalence of categories or a triequivalence of tricategories; see Figure~\ref{sdtwo}.

\item Comparing each degenerate bicategory with the monoidal category formed by its 1-cells and 2-cells, we exhibit an equivalence of categories of these structures, but not a biequivalence of bicategories or a triequivalence of tricategories; see Figure~\ref{sdthree}.

\end{itemize}

So to achieve an equivalence between the totalities of structures in question, we see that the ``correct" number of dimensions to take into account is critical.

Finally in Section~\ref{overall} we include some hypotheses about higher dimensions.

\enlargethispage{-6cm}

\subsection*{Acknowledgements}

We would like to thank Martin Hyland, Tom Leinster, Peter May, and Mike Shulman for their insightful suggestions and challenging comments.  In particular we would like to thank Tom Leinster for drawing our attention to the distinguished invertible elements in the first place.

\vspace{10cm}

\section{Methodology}\label{method}

In this section we outline the various ways in which we compare the structures in question.  We consider the general situation of comparing on the one hand $k$-fold degenerate $n$-categories and on the other hand $(n-k)$-categories with ``extra structure".  \emph{A priori} we have:

\begin{itemize}

\item an $(n+1)$-category \cat{nCat} of $n$-categories, $n$-functors, $n$-transformations, and so on

\item an $(n+1)$-category $\cat{nCat}(k)$ of $k$-fold degenerate $n$-categories, $n$-functors between them, $n$-transformations between those, and so on, as a full sub-$(n+1)$-category of \cat{nCat}

\item an $(n-k+1)$-category $\cat{PT}(n,k)$ of ``$k$-tuply monoidal $(n-k)$-categories", as hypothesised by the Periodic Table.

\end{itemize}

\noindent Our task is then to compare $\cat{nCat}(k)$ and $\cat{PT}(n,k)$; to do this, we regard $\cat{PT}(n,k)$ as a (partially discrete) $(n+1)$-category by adding in identity $j$-cells for all the ``missing" higher dimensions $n-k+2 \leq j \leq n+1$.  We can then look for an $(n+1)$-equivalence of $(n+1)$-categories
	\[\cat{nCat}(k) \map{\simeq} \cat{PT}(n,k).\]
(See Section~\ref{nequiv} for a definition of $(n+1)$-equivalence.)  

This approach does not produce a positive result in the cases studied.  On objects we ``forget" structure in the direction shown, but this does not necessarily give an $n$-functor.  So we examine lower-dimensional ``truncations" of the $(n+1)$-categories in question as follows.  For each $1 \leq j \leq n$ we write 

\begin{itemize}
\item $\cat{nCat}(k)_j$ for the $j$-dimensional ``truncation" of $\cat{nCat}(k)$
\item $\cat{PT}(n,k)_j$ to be the $j$-dimensional ``truncation" of $\cat{PT}(n,k)$
\end{itemize}
that is, the $j$-dimensional structures including only the lowest $j$-dimensions of the original $(n+1)$-categories.  Note that this is not, in general, a way of producing a $j$-category; we must check that $j$-cells compose strictly.  If so, we then determine whether the process of ``forgetting structure'' now gives a $j$-equivalence of $j$-categories.  We see that although truncation appears to be an unnatural process, this method produces positive results for a careful choice of $j$.

In fact, our very first task is to characterise \catk{nCat}{k}.  In detail, the various steps of the process are as follows.

\subsection{Precise characterisation of structures}\label{precise}

First we look at \catk{nCat}{k}\ which has
\begin{itemize}
\item 0-cells: $k$-degenerate $n$-categories i.e. $n$-categories with only one 0-cell, 1-cell, $\ldots$, $(k-1)$-cell
\item 1-cells: $n$-functors between these
\item 2-cells: $n$-transformations between those
\item 3-cells: $n$-modifications between those
\item 4-cells: $n$-perturbations between those
\item 5-cells: [no existing terminology]
\end{itemize}

\hspace{2em} $\vdots$

\begin{itemize}
\item $(n+1)$-cells: [no existing terminology]
\end{itemize}

\noindent We characterise each level of this structure precisely in terms of $(n-k)$-categories with extra structure, by taking the single hom-$(n-k)$-category and examining the remaining structure and axioms from the original degenerate $n$-category. 

We then consider that a $k$-degenerate $n$-category ``is precisely" an $(n-k)$-category with this extra structure; that is, one uniquely determines the other.  We continue in this fashion, producing statements of the form ``a $j$-cell in $\cat{nCat}(k)$ is precisely a $\ldots$" for all dimensions $j$.  





\subsection{Comparison functors}

The next stage is to look for comparison $j$-functors for each $j$-truncation
	\[\cat{nCat}(k)_j \lra \cat{PT}(n,k)_j.\]
In general a 0-cell on the left is one on the right together with some extra structure, so there is an obvious action in the direction shown, forgetting the extra structure.  We check that it is a $j$-functor, and ask if it is $j$-equivalence, using the recursive definitions of external and internal $j$-equivalence given below.  The forgetful action is canonical whereas the reverse direction involves \emph{choosing} extra structure.  

\label{nequiv}

\begin{Def} 
Let $x_1, x_2$ be 0-cells in a $j$-category $X$. 
\begin{itemize}
\item If $j=0$ then $x_1$ and $x_2$ are called \emph{internally equivalent} if and only if $x_1=x_2$.
\item If $j>0$ then $x_1$ and $x_2$ are called \emph{internally equivalent} (i.e. internally to the $j$-category $X$) if there are 1-cells
	\[x_1 \map{f} x_2 \mbox{\hspace{2em} \emph{and} \hspace{2em}} x_2 \map{g} x_1\]
such that $g \circ f$ is internally equivalent to $1_{x_1}$ in the hom-$(j-1)$-category $X(x_1,x_1)$ and $f \circ g$ is internally equivalent to $1_{x_2}$ in the hom-$(j-1)$-category $X(x_2,x_2)$.

\end{itemize}
\end{Def}

\begin{Def} \ \smallskip
\begin{itemize}
\item Let $X$ and $Y$ be 0-categories (i.e. sets). A 0-functor (i.e. function) $F:X \lra Y$ is called an \emph{external 0-equivalence} if and only if it is an isomorphism.
\item Let $j > 0$ and let $X$ and $Y$ be $j$-categories.  A $j$-functor $F:X \lra Y$ is called an \emph{external $j$-equivalence} or \emph{$j$-equivalence of $j$-categories} if
	\begin{enumerate}
	\item it is essentially surjective on 0-cells, i.e. given any 0-cell $y \in Y$ there is a 0-cell $x \in X$ such that $Fx$ is internally equivalent to $y$ in $Y$, and 
	\item it is locally a $(j-1)$-equivalence, i.e. given any 0-cells $x_1$ and $x_2$ in $X$, the $(j-1)$-functor
	\[X(x_1,x_2) \map{F} Y(Fx_1, Fx_2)\]
is a $(j-1)$-equivalence of $(j-1)$-categories.
	\end{enumerate}

\end{itemize}
\end{Def}

\noindent If we unravel the definitions we see that a $j$-functor $X \map{F} Y$ is an external $j$-equivalence if and only if

\begin{enumerate}

\item it is locally essentially surjective at all dimensions, i.e. essentially surjective on 0-cells and for $1 \leq m \leq j$ given any $m$-cell $\beta: Fx_1 \lra Fx_2 \in Y$ there is an $m$-cell $\alpha : x_1 \lra x_2 \in X$ such that $F\alpha$ is internally equivalent to $\beta$, and
\item it is locally faithful at the top dimension, i.e. for any pair of $j$-cells $\alpha_1, \alpha_2: x_1 \lra x_2 \in X$ 
	\[F\alpha_1 = F\alpha_2 \Rightarrow \alpha_1 = \alpha_2.\]
\end{enumerate}

\noindent It is useful to bear this ``unravelling'' in mind when showing that something is \emph{not} an equivalence.

Having found comparison functors in one direction we seek functors in the opposite direction
	\[\cat{PT}(n,k) \lra \cat{nCat}(k)\]
by making canonical choices of extra structure.

\subsection{Strictness and other ways of producing an equivalence}\label{methstrict}

A natural question to ask is: how can we alter $\cat{nCat}(k)$ to ``improve'' the situation?  One approach is to consider strict rather than weak structures, which elimates some of the problematic data.  Another more sutble approach is to restrict or alter the cells in $\cat{nCat}(k)$ with the specific goal of producing the required equivalence.  Neither approach gives a definitive solution; the first is unnecessarily strict and the second is in general intractable.

\subsection{Algebraic vs non-algebraic}

Since much of the problem arises from the structure constraints in the original $n$-category, we expect the situation to be better in ``non-algebraic" theories of $n$-category, where structure constraints are not actually specified.  Note that the Periodic Table was first described by Baez and Dolan, and the theory of $n$-categories proposed by these same authors \cite{bd1} is a non-algebraic theory.

Another non-algebraic theory is that of Street \cite{str2}.  The second named author has investigated the case of doubly degenerate bicategories and has checked that the following theorem holds.

\begin{thm}
The category of doubly degenerate weak Street 2-categories is equivalent to the category of commutative monoids.
\end{thm}

\noindent This uses the results of \cite{gur1}.

\subsection{Terminology}

\begin{itemize}

\item We generally use the adjectives ``strict", ``weak" and ``lax" to mean:
	\begin{center}
	\begin{tabular}{lcl}
	strict & -- & on the nose \\
	weak & -- & up to isomorphism \\
	lax & -- & up to non-invertible constraint cell
	\end{tabular}
	\end{center}

\item We generally use ``\ncat" to mean \emph{weak} \ncat, except as usual  weak 2-category is called a ``bicategory".  Thus instead of \catk{2Cat}{k} we have \catk{Bicat}{k}.  ``2-category" is usually reserved for the strict case except in the Periodic Table.  Similarly for tricategories and 3-categories.

\end{itemize}

\section{Degenerate categories}\label{dcs}

In this section we examine degenerate categories, that is, categories with only one object.  We show that these ``are precisely" monoids, the only non-canonical part of the correspondence being the choice of the single object.  To avoid this issue we will always pick our single object to be $\ast$.  We then examine the  full sub-2-category of the 2-category $\mathbf{Cat}$ whose 0-cells are these degenerate categories.

\subsection{Basic results}

The following result is well-known and consists of a routine rewriting of standard definitions.

\numarabic

\begin{thm}\label{dc}  \ \smallskip
\begin{enumerate}

\item A category \cl{C} with only one object $\ast$ is precisely a monoid $M_\cl{C}$ whose elements are the morphisms of \cl{C}: 

\begin{itemize}
\item multiplication in $M_\cl{C}$ is given by composition of morphisms in \cl{C}
\item the unit in $M_\cl{C}$ is given by the identity morphism in \cl{C}.
\end{itemize}
Associativity and unit axioms correspond to those for \cl{C}.

\item Extending the above correspondence, a functor $\cl{C} \map{F} \cl{D}$ is a precisely monoid homomorphism $M_\cl{C} \map{F} M_\cl{D}$.  Functoriality corresponds to preservation of the unit and multiplication for the monoid.

\item Extending the above correspondence, a natural transformation $\alpha: F \Rightarrow G$ is precisely a distinguished element $d_\alpha \in M_\cl{D}$ such that for all $x \in M_\cl{C}$

\begin{equation}\label{eq1}
d_{\alpha} \cdot Fx = Gx \cdot d_{\alpha}
\end{equation}

\begin{itemize}

\item The element $d_{\alpha}$ is the component of $\alpha$ at the single object of $\cl{C}$.

\item Equation \ref{eq1} corresponds to naturality of $\alpha$; the naturality square must commute for all morphisms in \cl{C} i.e. all elements of $M_\cl{C}$.

\end{itemize}

\end{enumerate}
\end{thm}

\subsection{Overall structure for degenerate categories}\label{dcoverall}

We now summarise the above results.  We see that a sensible result can be proved about the 1-dimensional structure formed by degenerate categories, but the corresponding result for the 2-dimensional structure fails.


We introduce the following notation for the totality of degenerate categories:

\begin{itemize}

\item Write $\cat{Cat}(1)_1$ for the full subcategory of \Cat\ whose objects are the degenerate categories.  

\item Write $\cat{Cat}(1)_2$ for the full sub-2-category of the 2-category $\mathbf{Cat}$ whose 0-cells are the degenerate categories.  

\end{itemize}

For the totality of monoids we use the following notation:

\begin{itemize}

\item Write \cat{Mnd} for the category of monoids and monoid homomorphisms

\item Write $\cat{Mnd}_2$ denote the (discrete) 2-category of monoids and monoid homomorphisms; the only 2-cells are identities.

\end{itemize}

\noindent Then there is a canonical functor
	\[\phi_1 : \catkj{Cat}{1}{1} \lra \cat{Mnd}\]
which ``forgets" the single object of each degenerate category.  However, there is no obvious canonical functor
	\[\catkj{Cat}{1}{2} \lra \cat{Mnd}_2\]
since this would have to send 2-cells with different source and target to identities on the right hand side; see remarks below.

\begin{thm} \label{dce}
$\phi_1$ gives an equivalence of categories.
\end{thm}

\begin{prf}
$\phi_1$ is clearly full, faithful and surjective on the nose.  Further, a (strict) inverse can be constructed by sending a given monoid $A$ to the corresponding degenerate category with single object $\ast$.  
\end{prf}

\subsubsection*{Remarks}

\numarabic

\begin{enumerate}
\item Note that a monoid $A$ can be realised as a degenerate category with \emph{any} one-element set as its set of objects.  If we do not choose to fix the single object to be $\ast$ then the fibre of $\phi_1$ over a given monoid $A$ is canonically isomorphic to the category of one-element sets, and we get a canonical pseudo-inverse to $\phi_1$ for each one-element set.  


\item Note that given the above inverse $\phi_1$, we can extend it to a (strict) 2-functor
	\[\cat{Mnd}_2 \lra \catkj{Cat}{1}{2}\]
by specifying its action on 2-cells.  
 However the resulting 2-functor is not locally full; it is straightforward to exhibit a transformation with a non-identity distinguished element.  

\end{enumerate}


\section{Doubly degenerate bicategories}\label{vdbs}

We now turn our attention to bicategories with only one 1-cell (and hence only one 0-cell); as in the previous section, we assume that the single 0-cell is $\ast$ and thus the single 1-cell is $I_\ast$.  We call these \emph{doubly degenerate bicategories}.

We begin by considering \cat{Bicat} with

\begin{itemize}
\item 0-cells: bicategories
\item 1-cells: weak functors
\item 2-cells: weak transformations
\item 3-cells: modifications
\end{itemize}
i.e. all structure constraints are invertible.  We then consider the full sub-tricategory of $\mathbf{Bicat}$ whose 0-cells are the doubly degenerate bicategories.   We also consider the lax (non-invertible) variant, which actually turns out to give the same results.  

\subsection{Basic results}

The following theorem is partly due to Leinster; part 1 is described in \cite{lei4} and part 2 in \cite{lei10}.

\numarabic

\begin{thm}\label{vdb}  \ \smallskip
\begin{enumerate}

\item A doubly degenerate bicategory \cl{B} with 0-cell $\ast$ is precisely a commutative monoid $X_\cl{B}$ equipped with a distinguished invertible element $d_{X} \in X$; we write $(X, d_{X})$ for this structure.

\item Extending the above correspondence, a weak functor $(X, d_{X}) \rightarrow (Y, d_{Y})$ is precisely a monoid homomorphism $F: X \rightarrow Y$ together with a distinguished invertible element $m_{F} \in Y$; we write $(F, m_{F})$.  Composition is given by $(G, m_{G}) \circ (F, m_{F}) = (GF, \ Gm_F \cdot m_{G})$.  Furthermore, all lax functors turn out to be weak.

\item Extending the above correspondence, a weak transformation 
	\[(F, m_{F}) \Rightarrow (G, m_{G})\] 
is the assertion that $F=G$ as monoid homomorphisms.  Furthermore, all lax transformations turn out to be weak.

\item A modification between such assertions then is precisely a distinguished element $\Gamma \in Y$ (which is not necessarily invertible). 

\end{enumerate}
\end{thm}

\noindent Thus doubly degenerate bicategories can be thought of as commutative monoids with some extra structure as above; we will eventually sum up these results in a theorem (analogous to Theorem~\ref{dce}) comparing these with ordinary commutative monoids without the extra structure.  We will exhibit a biequivalence at the 2-dimensional level but no equivalence at the 1- or 3-dimensional levels.  First we will provide the proofs of the four parts separately. 

\subsection{Bicategories}\label{bicatproof}

In this section we prove Theorem~\ref{vdb}, part 1, characterising doubly degenerate bicategories. By the results of Section~\ref{dcs} we know that the 2-cells form a monoid under $\circ$, and an Eckmann-Hilton type argument shows that this is commutative.  It is tempting to apply the Eckmann-Hilton argument to the operations $\circ$ and $\ast$ but some care must be taken to prove that $\ast$ is strictly unital; \emph{a priori} its unit acts only as strictly as the 1-cell units in the bicategory in question. 

Another approach is to define a new operation $\odot$, derived from $\ast$, which \emph{is} strictly unital, and apply the Eckmann-Hilton argument to $\odot$ and $\circ$.  In fact it is quite straightforward to prove directly that $\circ$ is commutative once one has considered the above argument using $\odot$, and this is the proof we give here.  

We then use the commutativity of $\circ$ to show that the operations $\circ$ and $\ast$ are the same, which enables us to show that the correspondence described in the theorem is indeed bijective.

Let \cl{B} be a doubly degenerate bicategory with only one 0-cell $*$ and only one 1-cell $I_*$. As usual we write $\curlyl, \curlyr, \curlya$ for the left and right unit and associativity constraints, omitting subscripts since there is only one 1-cell in any case. We construct from \cl{B} a commutative monoid and distinguished invertible element, written $(X, d_{X})$.  The underlying monoid is given by the degenerate hom-category $\cl{B}(*,*)$.


We now show that $\circ$ is commutative, using the following crucial facts:

\begin{enumerate}

\item $\curlyr = \curlyl$.  This is proved in \cite{js1}; alternatively it can be deduced from the coherence theorem for bicategories in the form ``all diagrams of constraints commute".

\item For any 2-cell $\alpha \in \cl{B}$ we have $\alpha = \curlyl \circ (1 \ast \alpha) \circ \curlyl^{-1}$ by naturality of $\curlyl$.  Similarly, $\alpha = \curlyr \circ ( \alpha \ast 1) \circ \curlyr^{-1}$.


\item The usual interchange law for $\ast$ and $\circ$
	\[(a \circ b) \ast (c \circ d) = (a \ast c) \circ (b \ast d).\]

\end{enumerate}

\noindent We then have the following calculation.

		\[\begin{array}{rcll}
	\beta \circ \alpha & = & (\curlyl \circ (1 \ast \beta) \circ \curlyl^{-1}) \circ (\curlyr \circ (\alpha \ast 1) \circ \curlyr^{-1}) & \mbox{by naturality of $\curlyl$ and $\curlyr$}\\
	&=& (\curlyl \circ (1 \ast \beta) \circ \curlyl^{-1}) \circ (\curlyl \circ (\alpha \ast 1) \circ \curlyl^{-1}) & \mbox{since $\curlyl=\curlyr$} \\
	&=& \curlyl \circ (1 \ast \beta) \circ (\alpha \ast 1) \circ \curlyl^{-1} & \mbox{since $\curlyl^{-1} \circ \curlyl = 1$} \\
	&=& \curlyl \circ ((1 \circ \alpha) \ast (\beta \circ 1)) \circ \curlyl^{-1} & \mbox{by interchange} \\
	&=& \curlyl \circ (\alpha \ast \beta) \circ \curlyl^{-1} & \mbox{[N.B. this is $\alpha \odot \beta$]}\\
 	&=& \curlyl \circ ((\alpha \circ 1) \ast (1 \circ \beta)) \circ \curlyl^{-1} \\
	&=& \curlyl \circ (\alpha \ast 1) \circ (1 \ast \beta) \circ \curlyl^{-1} & \mbox{by interchange} \\
	&=& (\curlyl \circ (\alpha \ast 1) \circ \curlyl^{-1}) \circ (\curlyl \circ (1 \ast \beta) \circ \curlyl^{-1}) & \mbox{since $\curlyl^{-1} \circ \curlyl = 1$}\\
	&=& (\curlyr \circ (\alpha \ast 1) \circ \curlyr^{-1}) \circ (\curlyl \circ (1 \ast \beta) \circ \curlyl^{-1}) & \mbox{since $\curlyr=\curlyl$} \\
	&=& \alpha \circ \beta & \mbox{by naturality of $\curlyl$ and $\curlyr$}
	\end{array}\]

\paragraph{N.B.} This calculation is essentially the Eckmann-Hilton process between $\circ$ and a new operation $\odot$ defined by 
	\[\alpha \odot \beta = \curlyl \circ (\alpha \ast \beta) \circ \curlyl^{-1}.\]

From the above calculation we see that 
	\[\alpha \circ \beta = \alpha \odot \beta\]
and that this operation commutes; thus we also have
	\[\begin{array}{rcll}
	\alpha \odot \beta &=& \curlyl \circ (\alpha \ast \beta) \circ \curlyl^{-1} \\
	&=& (\alpha \ast \beta) \circ \curlyl \circ \curlyl^{-1} & \mbox{by commutativity of $\circ$} \\
	&=& \alpha \ast \beta.
	\end{array}\]
This tells us that the two most ``obvious" multiplications we might define on the 2-cells of \cl{B}, by horizontal and vertical composition, are in fact the same and commutative.  Thus $X$ is a commutative monoid under $\circ$ (or $\ast$ or $\odot$) with unit given by $1_{I_*}$.

\paragraph{N.B.} It is worth noting that, \emph{a priori}, $\ast$ does not give a monoid structure on the 2-cells of $X$. It is tempting to argue that ``$\circ$ and $*$ give two monoid structures on the 2-cells of $X$ but by an Eckmann-Hilton argument these are the same and commutative".  This argument appears in the literature for \emph{strict} 2-categories, in which case $*$ \emph{does} give a monoid structure and the Eckmann-Hilton argument follows immediately.  

\bigskip
Continuing with the characterisation of doubly degenerate bicategories, we also have structure constraints
	\[\begin{array}{ccl}
	\curlya &:& (I \circ I) \circ I \stackrel{\sim}{\Rightarrow} I \circ (I \circ I),\\[3pt]
	\curlyl &:& I \circ I \stackrel{\sim}{\Rightarrow} I \mbox{\ \ and}\\[3pt]
	\curlyr &:& I \circ I \stackrel{\sim}{\Rightarrow} I.
	\end{array}\]
So \emph{a priori} this gives three distinguished invertible elements of $X$ corresponding to $\curlya, \curlyl$ and $\curlyr$.

However, we know that $\curlyr=\curlyl$; further, the pentagon identity gives us $\curlya^3 = \curlya^2$, and since $\curlya$ is invertible we have $\curlya=1$.  This leaves only one distinguished invertible element $\curlyl=\curlyr$ and we write $d_X = \curlyl = \curlyr$.  Thus we have a commutative monoid with a distinguished invertible element as asserted.

%


\subsection{Weak functors}\label{weakfunctors}

In this section we prove Theorem~\ref{vdb}, part 2, characterising weak functors between doubly degenerate bicategories.  We continue to use the results of Section~\ref{dcs} for degenerate categories, since all the hom-categories in the present situation are by definition degenerate.

A weak functor 
	\[(F,\phi): (X,d_X) \rightarrow (Y,d_Y)\] 
between doubly degenerate bicategories consists of the following.

\begin{itemize}


\item A functor on the unique hom-category; this is a degenerate category, so by Section~\ref{dcs} this gives us a monoid morphism
	\[ F : X \lra Y.\]

\item We have a structure constraint for composition.  Since there is only a single 1-cell, this gives a 2-cell isomorphism 
	\[\phi_{II}:FI \circ FI \stackrel{\sim}{\Rightarrow} F(I \circ I)\] 
in the target, satisfying naturality. This reduces to the condition 
	\[\phi_{II} \circ F(\beta \alpha) = F(\beta \alpha) \circ \phi_{II}\] 
for all $\alpha, \beta \in X$.  This automatically holds by commutativity, so the axiom gives us no further information.  Thus this data amounts to a distinguished invertible element $\phi_{II}$ in $Y$, which we call $m_{2}$.

\item There is a structure constraint for the unit.  This gives a 2-cell isomorphism 
	\[\phi_{*}: I'_{F*} \Rightarrow F(I_{*})\] 
in the target (subject to a vacuous naturality condition), that is another distinguished invertible element $m_{0} = \phi_{\ast} \in Y$.
\end{itemize}

\noindent Finally, we have three axioms for weak functors of bicategories.  The first is an associativity axiom which reduces to the equation ${m_2}^2 = {m_2}^2$.

There are two unit axioms, both giving
	\begin{equation}\label{eq4}
	d_{Y} = Fd_{X} \cdot m_{2} \cdot m_{0}.
	\end{equation}
We can rewrite (\ref{eq4}) as
	\begin{equation}\label{eq5}
	m_{0} = d_{Y} \cdot {m_2}^{-1} \cdot (Fd_X)^{-1}
	\end{equation}
so $m_{0}$ is determined by the rest of the data, leaving effectively just one distinguished invertible element that can be freely chosen.  We call this $m_{F} = m_2$ and we conclude that a weak functor gives a monoid homomorphism $F: X \rightarrow Y$ together with a distinguished invertible element $m_{F} \in Y$.


We now examine composition of weak functors.  Given functors
	\[X \map{(F,\phi^{F})} Y \map{(G,\phi^G)} Z\]
between doubly degenerate bicategories, with corresponding monoid maps and distinguished invertible elements
	\[(X,d_X) \map{(F,m_F)} (Y,d_Y) \map{(G,m_G)} (Z,d_Z),\]
it is routine to check that the composite corresponds to
	\[(G, m_{G}) \circ (F, m_{F}) = (GF,\ Gm_F \cdot m_{G}).\]  We observe that this composition is strictly associative and unital.

\paragraph*{Remark}\label{remlaxfun}

If we consider a lax functor instead of a weak one, then \textit{a priori} $m_{2}$ and $m_{0}$ are not invertible.  However, by equation (\ref{eq4}) we have 
	\[({d_Y}^{-1} \cdot Fd_X \cdot m_2) \cdot m_0 = 1\]
so by commutativity we have an inverse for $m_0$, and similarly for $m_2$.  Thus every lax functor between doubly degenerate bicategories is actually a weak functor and the situation is as above.   

\bigskip


\subsection{Weak transformations}

In this section we prove Theorem~\ref{vdb} part 3, characterising weak transformations.  So we consider a weak transformation $\sigma$ of doubly degenerate bicategories as shown below.

\[\xy
{\ar@/^1pc/^{(F,m_F)} (0,0)*+{(X,d_X)}; (20, 0)*+{(Y,d_Y)} };
{\ar@/_1pc/_{(G,m_G)} (0,0)*+{(X,d_X)}; (20,0)*+{(Y,d_Y)} };
(10,0)*{\Downarrow \sigma}
\endxy
\]
Such a weak transformation consists of an invertible 2-cell in the target which we write as $\sigma \in Y$.




The naturality condition for $\sigma$ gives
	\[F\alpha \cdot \sigma = \sigma \cdot G\alpha\]
for all $\alpha \in X$.  But we know that $Y$ is commutative and $\sigma$ invertible, so we have		
	\[F\alpha = G\alpha,\]  
i.e. $F=G$ as monoid homomorphisms.

Now we examine the two axioms for a weak transformation.  The first is the axiom for the associator which reduces to the equation 
	\[\sigma^{2}\cdot m_{F} = \sigma \cdot m_{G}\]
hence 
	\[\sigma \cdot m_{F} = m_{G}\]
since $\sigma$ is invertible.  So in fact $\sigma$ is completely determined by $m_{F}$ and $m_{G}$.

The axiom for units reduces to



	\[\sigma \cdot m_{F} \cdot Fd_{X} = m_{G} \cdot Gd_{X}.\] 
But we know $F=G$ and $\sigma = m_{G}\cdot {m_F}^{-1}$, so this equation is automatically satisfied and gives no further information.  

Hence a weak transformation $(F,m_F) \Rightarrow (G,m_G)$ is simply the assertion that $F=G$.

\paragraph{Remark}\label{remlaxtr}  Again, we note the result for the lax case.  From the second axiom we see that $\sigma$ is forced to be invertible, so every lax transformation of doubly degenerate bicategories is in fact a weak transformation, and the result is as above.

\subsection{Modifications}
We now examine modifications.  First, note that for a transformation 		
	\[(F,m_F) \Rightarrow (G,m_G)\]
to exist, we must have $F=G$, and in this case there is precisely one transformation, with ``component" $\sigma = m_{G} \cdot {m_F}^{-1}$.  So a modification must go from a transformation of this form to itself; it then consists of a 2-cell $\Gamma$ in the target bicategory such that
	\[\sigma \cdot \Gamma = \Gamma \cdot \sigma;\]
this holds for every $\Gamma$ since $Y$ is commutative.
So a modification ``between two assertions $F=G$" is simply a distinguished element $\Gamma \in Y$.

\subsection{Overall structure for doubly degenerate bicategories}
We will now summarise the above results and compare the above structures with ordinary commutative monoids (with no extra structure).  First, note that bicategories and weak functors form a category, but bicategories, weak functors, and weak transformations do not form a bicategory; adding in modifications does produce a tricategory.  The situation for doubly degenerate bicategories is better -- $\mathbf{Bicat}(2)_{2}$ is in fact a strict 2-category.

We write $\mathbf{Bicat}(2)_j$, for $j=1,2,3$, for the $j$-category of doubly degenerate bicategories consisting of (where appropriate):

\begin{itemize}
\item 0-cells: doubly degenerate bicategories
\item 1-cells: weak functors between them
\item 2-cells: weak transformations between those
\item 3-cells: modifications between those.
\end{itemize}

\noindent For the totality of commutative monoids we write 

\begin{tabular}{ll}
\cat{CMon} & for the category of commutative monoids and their homomorphisms \\
$\cat{CMon}_2$ & for the discrete 2-category on this category \\
$\cat{CMon}_3$ & for the discrete 3-category on this category.
\end{tabular}

\noindent Then we have for each $j=1,2,3$, a $j$-functor 
\begin{displaymath}
\xi_{j}:\mathbf{Bicat}(2)_j \longrightarrow \mathbf{CMon}_{j}
\end{displaymath}
and we have the following result.

\begin{thm}\label{vdbe} 
$\xi_2$ is a strict 2-equivalence of strict 2-categories, but $\xi_1$ and $\xi_3$ are not equivalences as they are not (locally) faithful.
\end{thm}

\begin{prf}
$\xi_2$ is evidently surjective on objects and locally surjective on 1-cells.  Moreover, it is locally an isomorphism on 2-cells, so it is a 2-equivalence.  $\xi_1$  is not faithful as it forgets the distinguished invertible element associated with a 1-cell in \catkj{Bicat}{2}{1}; similarly $\xi_3$ is not locally faithful at the 3-cell level as it forgets the distinguished element.  (As in Section~\ref{dcoverall} it is straightforward to exhibit a modification with a non-identity distinguished element.)
\end{prf}

\paragraph{Remark} A pseudo-inverse to $\xi_2$ is easily constructed by choosing distinguished invertible elements to be the identity.  

\bigskip

\subsection{Strictness}\label{vdbother}

Here are two possible restrictions on $\mathbf{Bicat}(2)_j$ in order to achieve the desired equivalences.

\begin{enumerate}
\item\label{option1} For $j=1$ onsider strict 2-categories and strict functors between them.  This gives an equivalence of categories.
\item In fact (\ref{option1}) is more restrictive than necessary; we can obtain an equivalence of categories by simply restricting the morphisms to those whose distinguished invertible element is the identity.  This amounts to considering those functors whose constraint $\phi_{II}$ is the identity; $\phi_{*}$ is then uniquely determined.
\end{enumerate}

However, no such methods are available to get a triequivalence for the tricategory $\mathbf{Bicat}(2)_{3}$; we would have to restrict to identity modifications.


\section{Degenerate bicategories}\label{dbs}

In this section we study degenerate bicategories, that is, bicategories with only one 0-cell $\ast$.  It is a well-known ``fact'' that monoidal categories ``are'' one-object bicategories; however, the bicategory of such things is more mysterious, as is the tricategory.

Of the following five parts in the following theorem parts 1 and 2 are well-known; part 3 was described by Leinster in \cite{lei10}.

\subsection{Basic results}\label{basic}

\begin{thm}\label{db} \ \smallskip

\begin{enumerate}

\item A bicategory with only one 0-cell $\ast$ is precisely a monoidal category.
\item Extending this correspondence, a weak functor $(F,\phi): X \rightarrow Y$ between such is precisely a weak monoidal functor.
\item A weak transformation $\alpha$ between such
\begin{displaymath}
\xy
{\ar@/^1pc/^{(F,\phi)} (0,0)*+{X}; (20, 0)*+{Y} };
{\ar@/_1pc/_{(G,\psi)} (0,0)*+{X}; (20,0)*+{Y} };
(10,0)*{\Downarrow \alpha}
\endxy
\end{displaymath}
is then precisely a distinguished object $\alpha$ in the monoidal category $Y$ together with, for every $A \in \textrm{ob\ }X$, an isomorphism
\begin{displaymath}
\alpha_{A}: GA \otimes \alpha \map{\sim} \alpha \otimes FA
\end{displaymath}
in $Y$, such that the following diagrams commute: first 
for any $A \map{f} B$ in $X$,   
\begin{displaymath}
\xymatrix{
GA \otimes \alpha \ar[r]^{\alpha_{A}} \ar[d]_{Gf \otimes 1} & \alpha \otimes FA \ar[d]^{1 \otimes Ff} \\
GB \otimes \alpha \ar[r]_{\alpha_{B}} & \alpha \otimes FB }
\end{displaymath}
secondly, 
\begin{displaymath}
\xymatrix{ 
(GA \otimes GB) \otimes \alpha \ar[r]^{\smallcurlya} \ar[ddd]_{\phi_{AB} \otimes 1} & GA \otimes (GB \otimes \alpha) \ar[r]^{1 \otimes \alpha_{B}} & GA \otimes (\alpha \otimes FB) \ar[r]^{a^{-1}} & (GA \otimes \alpha) \otimes FB \ar[d]^{\alpha_{A} \otimes 1} \\
&&& (\alpha \otimes FA) \otimes FB \ar[d]^{\smallcurlya} \\
&&& \alpha \otimes (FA \otimes FB) \ar[d]^{1 \otimes \phi_{AB}} \\
G(A \otimes B) \otimes \alpha \ar[rrr]_{\alpha_{A \otimes B}} &&& \alpha \otimes F(A \otimes B) }
\end{displaymath}
and finally
\begin{displaymath}
\xymatrix{ 
I \otimes \alpha \ar[r]^{\smallcurlyl} \ar[d]_{\phi \otimes 1} & \alpha \ar[r]^{\smallcurlyr^{-1}} & \alpha \otimes I \ar[d]^{1 \otimes \phi} \\
GI \otimes \alpha \ar[rr]_{\alpha_{I}} && \alpha \otimes FI }
\end{displaymath} A lax transformation is as above but without the requirement that $\alpha_{A}$ be invertible. 

\item A modification $\Gamma$ between such

\[
 \xy
{\ar@/^2pc/^{(F,\phi)} (0,0)*+{X}; (30,0)*+{Y}};
{\ar@/_2pc/_{(G,\psi)} (0,0)*+{X}; (30,0)*+{Y}};
{\ar@/_.5pc/@{=>}_{\alpha} (10,6)*+{}; (11,-6.5)*+{}};
{\ar@/^.5pc/@{=>}^{\beta} (20,6)*+{}; (19,-6.5)*+{}};
{\ar@3{->}^{\Gamma} (11,0)*{}; (19,0)*{} }
\endxy
\]
is then precisely a morphism $\Gamma: \alpha \rightarrow \beta$ in $Y$ (between distinguished elements) such that the following diagram commutes:
\begin{displaymath}
\xymatrix{
GA \otimes \alpha \ar[r]^{1 \otimes \Gamma} \ar[d]_{\alpha_{A}} & GA \otimes \beta \ar[d]^{\beta_{A}} \\
\alpha \otimes FA \ar[r]_{\Gamma \otimes 1} & \beta \otimes FA }
\end{displaymath}


\end{enumerate}
\end{thm}

\begin{prf}
The proof of the above results consists of a routine rewriting of the definitions.
\end{prf}

\subsection{Overall structure for degenerate bicategories}\label{dboverall}

As in the previous sections, we now compare the above structures with ordinary monoidal categories; we will exhibit an equivalence at the level of categories but not at any of the higher-dimensional possibilities.


More precisely, for $j=1,2,3$ write $\cat{Bicat}(1)_j$ for the $j$-dimensional structure consisting of (where appropriate)
\begin{itemize}
\item 0-cells: degenerate bicategories
\item 1-cells: weak functors between them
\item 2-cells: weak transformations between those
\item 3-cells: modifications between those.
\end{itemize}

Then for $j=1$ we have a category and for $j=3$ a tricategory. 
(We do not explicitly give the tricategorical structure, regarding it as inherited from the tricategory $\mathbf{Bicat}$.)  However, for $j=2$ we do not have a bicategory as composition of 2-cells is not strictly associative or unital. 

For the totality of monoidal categories we write $\mathbf{MonCat}_j$ for the $j$-category consisting of (where appropriate)
\begin{itemize}
\item 0-cells: monoidal categories
\item 1-cells: monoidal functors
\item 2-cells: monoidal tranformations
\item 3-cells: identities.
\end{itemize}

\begin{thm}
The obvious functor
	\[\xi : \cat{Bicat}(1)_1 \rightarrow \mathbf{MonCat}_1\]
is an equivalence.  
\end{thm}

\subsection{Other ways of producing an equivalence for degenerate bicategories}\label{dbother}

As in the doubly degenerate case, one natural question to ask is:  how can we alter $\mathbf{Bicat}(1)$ in order to ``improve'' the situation?  Recall that


\begin{itemize}

\item in $\cat{Bicat}(1)_3$ a 2-cell $F \lra G$ has components which are isomorphisms
	\[GA \otimes \alpha \map{\sim} \alpha \otimes FA\] whereas 
\item in $\cat{MonCat}_3$ a 2-cell $F \lra G$ has components
	\[FA \lra GA.\]
\end{itemize}

\noindent Thus in order to get a comparison functor at all, we consider oplax rather than weak transformations.  We now discuss two possible strategies.
\begin{enumerate}
\item Restrict the 2-cells in $\mathbf{Bicat}(1)_{3}$ to those having their distinguished element the unit.  However, this is closed under neither vertical nor horizontal composition.  In principle we could try to take the closure under composition; in practice this is both difficult and unilluminating.
\item We could consider the ``backwards'' direction and at least get a functor
\[
\mathbf{MonCat}_{3} \rightarrow \mathbf{Bicat}(1)_{3}
\]
as follows.  We send a 2-cell $\alpha: F \Rightarrow G$ to the transformation
with distinguished element $I$ and components $\alpha_{A}: I \otimes FA \rightarrow GA \otimes I$ given by
\begin{equation}\label{switching}
\xymatrix{
I \otimes FA \ar[r]^{\smallcurlyl} & FA \ar[r]^{\alpha_{A}} & GA \ar[r]^{\smallcurlyr^{-1}} & GA \otimes I. }
\end{equation}
To show that this is not locally essentially surjective on 2-cells, it is once again straightforward to exhibit a transformation whose distinguished element is not isomorphic to $I$.
\end{enumerate}

\subsection{Example: monad functors}

In this section we discuss the example of monad functors.  This example makes the transformations between degenerate bicategories look a little less mysterious.  We will see that monads, their functors, and their transformations are examples of functors, transformations, and modifications between degenerate bicategories.  The idea is that the bicategory of monads in an arbitrary bicategory $\cl{B}$ on a fixed object $X$ is given by the hom-bicategory $\mathbf{Bicat}_{l}(1, \cl{B}_{X})$, where the subscript $l$ indicates that we take lax functors and lax transformations, and $\cl{B}_{X}$ is the full sub-bicategory of $\cl{B}$ with the single object $X$.  Here, as usual, 1 indicates the terminal bicategory.  The bicategories 1 and $\cl{B}_{X}$ are manifestly degenerate.




Unpacking this presentation, we get the usual definitions as follows.  Let
	\[\begin{array}{ccccc}
	S & : & \cl{C} & \lra & \cl{C} \\
	T & : & \cl{D} & \lra & \cl{D}
	\end{array}\]
be monads.  Then

\begin{itemize}
\item a monad functor $S \lra T$ consists of a functor
	\[U : \cl{C} \lra \cl{D}\]
together with a natural transformation
	\[\phi : TU \lra US\]
satisfying two axioms.  Here $U$ plays the role of the distinguished element.

\item A monad functor transformation 
	\[(U,\phi) \lra (U', \phi')\]
is then a natural transformation 
	\[\Gamma : U \lra U'\]
such that the following diagram commutes.

\[\xy
{\ar^{T\Gamma} (0,0)*+{TU}; (20,0)*+{TU'} };
{\ar_{\phi} (0,0)*+{TU}; (0,-14)*+{US} };
{\ar^{\phi'} (20,0)*+{TU'}; (20,-14)*+{U'S} };
{\ar_{\Gamma S} (0,-14)*+{US}; (20,-14)*+{U'S} };
\endxy
\]
Thus $\Gamma$ arises as a morphism of distinguished elements.
\end{itemize}

\subsection{Example: topological analogue}

The results of the previous sections indicate that the top-dimensional cells in the $(n+1)$-category $\cat{nCat}(k)$ cause unavoidable ``problems" when it comes to looking for equivalences with the structures given by the Periodic Table.  In this section we discuss a topological analogue to suggest that this is not ``wrong" but a phenomenon that does arise naturally elsewhere.  Topology provides a natural example in which the hom-$n$-category $\mathbf{nCat}(k)(X,Y)$ between two $k$-degenerate $n$-categories has interesting top-dimensional cells.  We indicate why the topological analogues of doubly degenerate bicategories should form a tricategory with nontrivial 3-cells, by computing the homotopy groups of their mapping space.

From one point of view, weak $\w$-groupoids ``are'' spaces and weak $n$-groupoids ``are'' $n$-truncated spaces, i.e., spaces $X$ with $\pi_{r}X=0$ if $r>n$.  Using this analogy, the weak $\w$-groupoid of functors $X \rightarrow Y$ between weak $\w$-groupoids should correspond to the (unbased) mapping space functor $\textrm{Map}(X,Y)$ for spaces.  This is a viable position for $n$-groupoids as well, as shown by the following proposition.

\begin{prop}
Let $Y$ be a connected, $n$-truncated space, and $X$ any CW-complex.  Then $\textrm{Map}(X,Y)$ is $n$-truncated as well.
\end{prop}
\begin{proof}
Since $X$ is a CW-complex, there is a cofibration $X^{q} \hookrightarrow X^{q+1}$, where $X^{q}$ denotes the $q$-skeleton of $X$.  This induces a fibration
\begin{displaymath}
\textrm{Map}(X^{q+1}/X^{q}, Y) \rightarrow \textrm{Map}(X^{q+1},Y) \rightarrow \textrm{Map}(X^{q},Y).
\end{displaymath}
Using this fibration, we will prove by induction that $\pi_{r} \textrm{Map}(X^{q},Y) = 0$ when $r>n$, and thus that $\pi_{r} \textrm{Map}(X,Y)=0$ for all $r>n$.

When $q=0$, $\textrm{Map}(X^{0},Y) = \textrm{Map}(\coprod *,Y) \cong \prod Y$.  Since $\pi_{r}Y = 0$ when $r>n$, the same holds for the product.

Assume the result for $q$ by induction.  Then the fibration above gives a long exact sequence of homotopy groups, part of which is displayed below.

\begin{displaymath}
\pi_{r+1} \textrm{Map}(X^{q},Y) \rightarrow \pi_{r} \big( \textrm{Map}(X^{q+1}/X^{q}, Y) \big) \rightarrow \pi_{r} \textrm{Map}(X^{q+1},Y) \rightarrow \pi_{r} \textrm{Map}(X^{q},Y)
\end{displaymath}

\noindent By induction, the first and last groups are zero, giving an isomorphism of the middle two groups.  Since $X^{q+1}/X^{q}$ is a wedge of $(q+1)$-spheres, we have reduced the problem to computing $\pi_{r}\textrm{Map}(S^{q+1}, Y)$ for $r>n$; if this group is zero then we are finished.

There is another fibration
\begin{displaymath}
\textrm{Map}_{*}(S^{q+1}, Y) \rightarrow \textrm{Map}(S^{q+1}, Y) \rightarrow Y,
\end{displaymath}

\noindent where $\textrm{Map}_{*}$ is the space of based maps and the map $\textrm{Map}(S^{q+1}, Y) \rightarrow Y$ is the map induced by evaluation at the basepoint.  The long exact sequence in homotopy groups gives
\begin{displaymath}
\xymatrix{
\pi_{r+1}Y \ar[r] & \pi_{r} \textrm{Map}_{*}(S^{q+1},Y) \ar[r] & \pi_{r} \textrm{Map}(S^{q+1},Y) \ar[r] & \pi_{r}Y, }
\end{displaymath}
and the first and last groups in this sequence are zero since $r>n$ and $Y$ is $n$-truncated.  Thus the middle groups are isomorphic, and we compute that $\pi_{r} \textrm{Map}_{*}(S^{q+1},Y)$ is
\begin{displaymath}
[S^{r}, \textrm{Map}_{*}(S^{q+1}, Y)] \cong [S^{r+q+1}, Y] = 0
\end{displaymath}
by adjunction and the fact that $Y$ is $n$-truncated.
\end{proof}

\begin{cor}  If $X,Y$ are connected, $n$-truncated CW-complexes, then the space $\textrm{Map}(X,Y)$ is $n$-truncated as well.
\end{cor}

We can now see how the topology of $n$-truncated spaces predicts the non-triequivalence of the tricategory $\mathbf{Bicat}(1)_3$ and the tricategory $\mathbf{CMon}_{3}$.  First, we must restrict to groupoids; thus commutative monoids become abelian groups.  Since we are interested in doubly degenerate bigroupoids, we take as topological analogues the Eilenberg-Mac Lane spaces $K(G, 2)$'s.  These spaces have homotopy groups $\pi_{i}K(G,2) = 0$ if $i \neq 2$ and $\pi_{2} K(G,2) = G$; this single homotopy group characterizes an Eilenberg-Mac Lane space up to homotopy equivalence if we assume them to be CW-complexes.

If $A,B$ are abelian groups, then the hom-bigroupoid $\mathbf{CMon}_{3}(A,B)$ has only identity 2-cells.  Thus the space associated with this bigroupoid would have vanishing $\pi_{2}$.  We will show that this is not the case for the mapping space $\textrm{Map} \big( K(A,2), K(B,2) \big)$, thus the topology has predicted that $\mathbf{CMon}_{3}$ is the ``wrong'' tricategory of commutative monoids.

Consider the fibration 
\begin{displaymath}
\textrm{Map}_{*} \big( K(A,2), K(B,2) \big) \rightarrow \textrm{Map} \big( K(A,2), K(B,2) \big) \rightarrow K(B,2)
\end{displaymath}
as in the Proposition, where the second map is induced by the inclusion of the basepoint $* \hookrightarrow K(A,2)$.  The first term is homotopy equivalent to the discrete space of group homomorphisms $A \rightarrow B$.  Therefore the long exact sequence of homotopy groups induces an isomorphism $\pi_{i}\textrm{Map} \big( K(A,2),K(B,2) \big) \cong \pi_{i}K(B,2)$ for $i > 1$.  When $i=2$, we see that
\begin{displaymath}
\pi_{2} \textrm{Map} \big( K(A,2), K(B,2) \big) \cong B;
\end{displaymath}
\noindent this is the space-level analogue of our result that the 2-cells of the hom-bicategory $\mathbf{Bicat}(1)_3(X,Y)$ correspond to elements of $Y$.

%

\section{Higher-dimensional hypotheses}\label{overall}

In this section we further consider the question: how many dimensions of structure give us the equivalence we seek to make the Periodic Table precise?  Recall the answers obtained in the previous sections:

\begin{itemize}
\item degenerate categories: 1-dimensional structure instead of a possible 2
\item degenerate bicategories: 1-dimensional structure instead of a possible 3
\item doubly degenerate bicategories: 2-dimensional structure instead of a possible 3.
\end{itemize}

It would be desirable to give an answer for every entry of the Periodic Table, i.e., to be able to give a general answer for $k$-degenerate $n$-categories.  While this is currently far beyond our scope, we will make a small hypothesis in this direction.

First we observe that, of the above three cases, only the third (doubly degenerate bicategories) is a ``stable" case.  We suspect that this makes a difference when answering the above question, and in particular that in the stable cases the situation as a whole is more tractable.

Our hypothesis concerns the stable cases in the first column of the Periodic Table.  Recall that this column reads as follows:

\bigskip
\begin{tabular}{lcll}
1-degenerate categories & $\cong$ & monoids \\
2-degenerate 2-categories & $\cong$ & commutative monoids & \hspace{2em} {\sc stable}\\
3-degenerate 3-categories & $\cong$ & commutative monoids & \hspace{3em} $\downarrow$\\
4-degenerate 4-categories & $\cong$ & commutative monoids \\
& $\vdots$ \\
$n$-degenerate $n$-categories & $\cong$ & commutative monoids \\
& $\vdots$ 
\end{tabular}

\bigskip
Our hypothesis concerns the $(n+1)$-category \cat{nCat($n$)} and its truncations, with

\begin{itemize}
\item 0-cells: $n$-degenerate $n$-categories i.e. $n$-categories with only one $(n-1)$-cell
\item 1-cells: $n$-functors between these
\item 2-cells: $n$-transformations between those
\item 3-cells: $n$-modifications between those
\item 4-cells: $n$-perturbations between those
\item 5-cells: [no existing terminology]
\end{itemize}

\hspace{2em} $\vdots$

\begin{itemize}
\item $(n+1)$-cells: [no existing terminology]
\end{itemize}

The following three hypotheses extend what we have already proved for $n=2$.

\bigskip

\begin{hyp} {\bfseries \ Basic results}

\bigskip
\begin{em}
\noindent Let $n \geq 3$.

\begin{itemize}

\item 0-cells: An $n$-degenerate $n$-category is precisely a commutative monoid with $D(n)$ distinguished invertible elements.  We expect $D(n)$ to be a sequence of natural numbers that increase rapidly as $n$ increases, since the distinguished invertible elements arise as constraint $n$-cells in the original $n$-category.

\item 1-cells: An $n$-functor $f: X \lra Y$ between such is precisely a monoid homomorphism $f: X \lra Y$ with $F(n)$ distinguished invertible elements in $Y$.  We expect $F(n)$ also to be an increasing sequence.

\item 2-cells: An $n$-transformation $\alpha: f \Rightarrow g$ is the assertion $f=g$ as monoid homomorphisms, with no condition on distinguished invertible elements.  

\item 3-cells: An $n$-modification between such is the identity.

\item 4-cells: An $n$-perturbation between such is the identity.\\

\vspace{-2ex}\hspace{2em} $\vdots$

\item $(n+1)$-cells: an $(n+1)$-cell between such is a distinguished element in $Y$, not necessarily invertible.  This is because at this level the data will be ``for every 0-cell $x \in X$ an $n$-cell $\sigma_x \in Y$".  Since there is only one 0-cell $*$, we simply have one distinguished element in $Y$.  It will be required to satisfy some equations which we expect to give no further information since multiplication in $Y$ is commutative.

\end{itemize}
\end{em}
\end{hyp}

As in the previous sections, we sum up these results by considering the overall structure.  We write \cat{CMnd} for the category of commutative monoids and monoid homomorphisms.  For $j \geq 2$ we write $\cat{CMnd}_j$ for this category regarded as a discrete $j$-category by adding higher identity cells, and $\cat{CMnd}_1 = \cat{CMnd}$.

We write \cat{nCat($n$)} for the $(n+1)$-category of $n$-degenerate $n$-categories.  We write $\cat{nCat($n$)}_j$ for the $j$-truncation of this $(n+1)$-category, and $\cat{nCat($n$)}_{n+1} = \cat{nCat($n$)}$.  

\bigskip

\begin{hyp} {\bfseries \ Overall structure}

\bigskip
\noindent \emph{Let $n \geq 3$.  For $1 \leq j \leq n+1$ there is a forgetful $j$-functor
	\[\cat{nCat($n$)}_j \lra \cat{CMnd}_j.\]
This is not a $j$-equivalence for $j=1,n+1$, but is a $j$-equivalence for $2 \leq j \leq n$.  }
\end{hyp}


\bigskip

Finally we consider the question of eliminating the distinguished invertible elements by using a stricter form of $n$-category.  Generalising from the previous sections, we see that we do not need to restrict all the way to strict \ncats\ -- a semistrict version will suffice.  One form of semistrictness has everything strict except interchange (cf. $\mathbf{Gray}$-categories and see \cite{cra1,cra2}); another has everything strict except units \cite{koc1,sim3}.  These have both been proposed as solutions to the coherence problem for $n$-categories.


However, there are other possible ``shades" of semistrictness and the above notions do not appear to be right for the present purposes.  Instead, we need a form of semistrict \ncat\ in which the units and interchange for $(n-1)$-cells are strict, but everything else can be weak.  This is to eliminate the constraint $n$-cells that become distinguished invertible elements in our $n$-degenerate situation; we expect that as in the case $n=2$ the associator is automatically forced to be the identity.

\bigskip
\begin{hyp} {\bfseries \ Semistrictness}

\bigskip
\noindent \emph{Let $n \geq 3$.  Then an $n$-degenerate semistrict $n$-category in the above sense is precisely a commutative monoid.}

\end{hyp}

%

\providecommand{\bysame}{\leavevmode\hbox to3em{\hrulefill}\thinspace}
\providecommand{\MR}{\relax\ifhmode\unskip\space\fi MR }
\providecommand{\MRhref}[2]{%
  \href{http://www.ams.org/mathscinet-getitem?mr=#1}{#2}
}
\providecommand{\href}[2]{#2}

\end{document}